\documentstyle[12pt]{amsart}
\newtheorem{theorem}{Theorem}\newtheorem{lemma}{Lemma}
\newtheorem{fact}{Fact}
\newtheorem{definition}{Definition}

\renewcommand{\text}[1]{\quad\mbox{#1}\quad}
\let\pont\relax
\def\eps{\varepsilon}
\def\qed{\ifhmode\unskip\nobreak\fi\quad\ifmmode\Box\else$\Box$\fi}

\title{On $k$-ordered Hamiltonian Graphs}

\author{Gabor N. Sarkozy}
\thanks{This paper was written while 
Sarkozy was visiting
MSRI Berkeley, as part of the 
Combinatorics Program. Research at MSRI is supported
in part by NSF grant DMS-9022140.} 
\address{\hskip-\parindent Gabor N. Sarkozy\\
Computer Science Department\\Worcester
Polytechnic Institute\\Worcester, MA 01609}
\email{gsarkozy@@cs.wpi.edu}
\author{Stanley Selkow}
\address{\hskip-\parindent Stanley Selkow\\
Computer Science Department\\Worcester
Polytechnic Institute\\Worcester, MA 01609}
\email{sms@@cs.wpi.edu}

\begin{document}

\begin{abstract}
A Hamiltonian graph $G$ of order $n$ is $k$-ordered, $2\leq k \leq n$,
if for every sequence $v_1, v_2, \ldots ,v_k$ of $k$ distinct vertices
of $G$, there exists a Hamiltonian cycle that encounters $v_1, v_2, 
\ldots , v_k$ in this order. In this paper, answering a question
of Ng and Schultz, we give a sharp bound for the minimum degree 
guaranteeing that a graph is a
$k$-ordered Hamiltonian graph under some mild restrictions. 
More precisely, we show that there
are $\eps, n_0> 0$ such that if $G$ is a graph of order $n\geq n_0$
with minimum degree at least $\lceil \frac{n}{2} \rceil
+ \lfloor \frac{k}{2} \rfloor  - 1$ and
$2\leq k \leq \eps n$, then $G$ is a $k$-ordered Hamiltonian graph.
It is also shown that this bound is sharp for every $2\leq k \leq 
\lfloor \frac{n}{2} \rfloor$.
\end{abstract}
\maketitle

\section{Introduction}
\subsection{Notations and definitions}
For basic graph concepts see the monograph of Bollob\'as \cite{B}.\\
$+$ will sometimes be used for disjoint union of sets.
$V(G)$ and $E(G)$ denote the vertex-set and the edge-set
of the graph $G$. $(A,B,E)$ denotes a bipartite graph $G=(V,E)$,
where $V=A+B$, and $E\subset A\times B$.
For a graph $G$ and a subset $U$ of its vertices,
$G|_U$ is the restriction to $U$ of $G$.
$N(v)$ is the set of neighbours of $v\in V$.
Hence the size of $N(v)$ is $|N(v)|=deg(v)=deg_G(v)$,
the degree of $v$.
$\delta(G)$ stands for the minimum, and
$\Delta(G)$ for the maximum  degree in $G$.
$\nu(G)$ is the size of a maximum matching in $G$.
For $A\subset V(G)$ we write
$N(A)=\cap_{v\in A}N(v)$, the set of common neighbours.
$N(x,y,z,...)$ is shorthand for $N(\{x,y,z,...\})$.
For a vertex $v\in V$ and set $U\subset V-\{v\}$, we write
$deg(v,U)$ for the number of edges from $v$ to $U$.
When $A,B$ are disjoint subsets of $V(G)$,
we denote by $e(A,B)$ the number of edges of $G$
with one endpoint in $A$ and the other in $B$.
For non-empty $A$ and $B$,$$d(A,B)=\frac{e(A,B)}{|A||B|}$$
is the {\bf density} of the graph between $A$ and $B$.
\begin{definition}\pont
The bipartite graph $G=(A,B,E)$ is $\eps$-{\bf regular} if
$$X\subset A,\ Y\subset B,\ |X|>\eps|A|,\ |Y|>\eps|B|
\text{imply}|d(X,Y)-d(A,B)|<\eps,$$
otherwise it is $\eps$-irregular.
\end{definition}
We will often say simply that ``the pair $(A,B)$ is $\eps$-regular''
with the graph $G$ implicit.
We will also need a stronger version.
\begin{definition}\pont $(A,B)$ is $(\eps,\delta)$
{\bf super-regular} if it is $\eps$-regular and
$$deg(a)>\delta|B|\ \forall a\in A,\quad
deg(b)>\delta|A|\ \forall b\in B.$$
\end{definition}

\subsection{$k$-ordered Hamiltonian graphs}

Let $G$ be a graph on $n\geq 3$ vertices. A {\bf Hamiltonian cycle}
({\bf path})
of $G$ is a cycle (path) 
containing every vertex of $G$. A {\bf Hamiltonian graph}
is a graph containing a Hamiltonian cycle. A classical result of Dirac
\cite{D} asserts that if $\delta(G)\geq n/2$, then $G$ is Hamiltonian.
A {\bf Hamiltonian-connected} graph is a graph in which every pair of
vertices can be connected with a Hamiltonian path.

The following interesting 
concept was created by Chartrand: For a positive integer $2\leq k \leq n$,
and for a sequence $S = v_1, v_2, 
\ldots ,v_k$ of $k$ distinct vertices, a cycle $C$ in $G$ is called a
$v_1-v_2- \ldots -v_k$-cycle, or shortly an $S$-cycle, if 
the vertices of $S$ are encountered on $C$ in the
specified order. For a Hamiltonian graph $G$
we say that $G$ is {\bf $k$-ordered} if for every sequence $S = v_1, v_2, 
\ldots ,v_k$ of $k$ distinct vertices, there exists a Hamiltonian $S$-cycle.
It is not hard to see that every Hamiltonian
graph is both 2-ordered and 3-ordered. Furthermore, a Hamiltonian
graph $G$ of order $n$ is $n$-ordered if and only if $G=K_n$. Also, if
$G$ is $k$-ordered, then $G$ is $l$-ordered for every $2\leq l \leq k$
(see \cite{NS}).

A natural question is whether we can obtain a Dirac type condition on
the minimum degree for guaranteeing that the graph is a $k$-ordered
Hamiltonian graph. Indeed, the first result of this type was obtained in
\cite{NS}. In this paper it was shown (among other results) that if
$3\leq k \leq n$ and $\delta(G) \geq \frac{n}{2} + k - 3$, 
then $G$ is a $k$-ordered Hamiltonian
graph. The authors raised the question whether this can be improved.
In this paper our goal is to determine the best possible bound under
the restrictions that $n$ is sufficiently large and $k$ does not grow too
quickly with $n$. More precisely, our main result is the following.

\begin{theorem}\pont\label{tetel1}
There exist $\kappa, n_0> 0$
such that if a graph $G$ has order $n\geq n_0$ with
\begin{equation}\label{fok}
\delta(G)\geq \lceil \frac{n}{2} \rceil
+ \lfloor \frac{k}{2} \rfloor  - 1  ,\end{equation}
and $2\leq k \leq \kappa n$,
then $G$ is a $k$-ordered Hamiltonian graph.
\end{theorem}

It is not hard to see that this bound is best possible for most $k$-s.

\begin{theorem}\pont\label{tetel2}
For any $2\leq k \leq \lfloor \frac{n}{2} \rfloor$
there exists a graph $G$ of order $n$ with
\begin{equation}
\delta(G) = \lceil \frac{n}{2} \rceil
+ \lfloor \frac{k}{2} \rfloor  - 2,\end{equation}
such that $G$ is not a $k$-ordered Hamiltonian graph.
\end{theorem}

\section{The main tools}
In the proof the following lemma of Szemer\'{e}di
plays a central role.
\begin{lemma}[Regularity Lemma \cite{Sz}]\pont\label{reg}
For every positive $\eps$ and positive integer $m$
there are positive integers $M$ and $n_1$ with the following property:
for every graph $G$ with $n\geq n_1$ vertices
there is a partition of the vertex set
into $l+1$ classes (clusters)$$V=V_0+V_1+V_2+...+V_l$$such that
\begin{itemize}
\item\quad $m\leq l\leq M$
\item\quad $|V_1|=|V_2|=...=|V_l|$
\item\quad $|V_0|<\eps n$
\item\quad at most $\eps l^2$ of the pairs $\{V_i,V_j\}$
are $\eps$-irregular.
\end{itemize}
\end{lemma}
We will use the following simple consequence of Lemma \ref{reg}.

\begin{lemma}[Degree form]\pont\label{degree}
For every $\eps>0$ there is an $M=M(\eps)$
such that if $G=(V,E)$ is any graph
and $d\in[0,1]$ is any real number,
then there is a partition of the vertex-set $V$
into $l+1$ clusters $V_0,V_1,...,V_l$,
and there is a subgraph $G'=(V,E')$
with the following properties:
\begin{itemize}
\item $l\leq M$,
\item $|V_0|\leq\eps|V|$,
\item all clusters $V_i,\,i\geq1,$ are of the same size
$L\leq\lceil\eps|V|\rceil$.
\item $deg_{G'}(v)>deg_G(v)-(d+\eps)|V|\text{for all}v\in V$,
\item $G'|_{V_i}=\emptyset$ ($V_i$ are independent in $G'$),
\item all pairs $G'|_{V_i\times V_j},\ 1\leq i<j\leq l$,
are $\eps$-regular, each with a density either 0 or exceeding $d$.
\end{itemize}
\end{lemma}
The other main tool asserts that if $(A, B)$ is a super-regular
pair with $|A|=|B|$ and $x\in A, y\in B$, then there is a Hamiltonian path
starting with $x$ and ending with $y$. This is a very special case
of the Blow-up Lemma \cite{KSSz-BL}, and it has already appeared
in \cite{KSSz-B} and in \cite{S}. More precisely.
\begin{lemma}\pont\label{ham}
For every $\delta>0$ there are $\eps_0, n_2>0$ such that if
$\eps \leq \eps_0$ and $n\geq n_2$, $G=(A,B)$ is an $(\eps, \delta)$ 
super-regular pair with $|A|=|B|=n$ and $x\in A$, $y\in B$, then there is a 
Hamiltonian path in $G$ starting with $x$ and ending with $y$.
\end{lemma}

We will also use two simple P\'{o}sa-type lemmas on Hamiltonian-connectedness.
The second one is the bipartite version of the first one.
\begin{lemma}[see \cite{BE}]\pont\label{hamcon}
Let $G$ be a graph on $n\geq 3$ vertices with degrees
$d_1\leq d_2 \leq \ldots \leq d_n$ such that for every
$2\leq k \leq \frac{n}{2}$ we have $d_{k-1} > k$. Then $G$ is
Hamiltonian-connected.
\end{lemma}

\begin{lemma}[see \cite{BE}]\pont\label{bihamcon}
Let $G=(A,B)$ be a bipartite graph 
with $|A|=|B|=n\geq 2$ with degrees
$d_1\leq d_2 \leq \ldots \leq d_n$ from $A$ and with degrees
$d_1'\leq d_2' \leq \ldots \leq d_n'$ from $B$. Suppose that for every
$2\leq j \leq \frac{n+1}{2}$ we have $d_{j-1} > j$
and that for every
$2\leq k \leq \frac{n+1}{2}$ we have $d_{k-1}' > k$. Then $G$ is
Hamiltonian-connected.
\end{lemma}

Finally we will use the following simple fact.
\begin{lemma}[Erd\H{o}s, P\'{o}sa, see \cite{B}]\pont\label{ep}
Let $G$ be a graph on $n$ vertices. Then
$$\nu(G) \geq \min \{ \delta(G), \frac{n-1}{2} \}.$$
\end{lemma}

In case we have a good upper bound on the maximum degree of $G$, we
can strengthen this lemma in the following way.
\begin{lemma}\pont\label{matching}
In a graph $G$ of order $n$
$$\nu(G)\geq\delta(G)\frac{n}{2(\delta(G)+\Delta(G))}\geq
\delta(G)\frac{n}{4\Delta(G)}.$$
\end{lemma}

In fact, let us take a maximal matching $M$ with $m$ edges.
Then for the number of edges $E$ between $M$ and
$V(G)\setminus M$ we get
$\delta(G)(n-2m)\leq E\leq2m\Delta(G)$,
which proves the lemma.

\section{Proof of Theorem 1}
\subsection{Outline of the proof}

During the past couple of years the first author, together
with J. Koml\'{o}s and E. Szemer\'{e}di, developed a new
method
in graph theory based on the Regularity Lemma and the Blow-up Lemma. 
The method is usually
applied to find certain spanning subgraphs in dense graphs. Typical examples
are spanning trees 
(Bollob\'{a}s-conjecture, see \cite{KSSz-B}),
Hamiltonian cycles or powers of Hamiltonian cycles
(P\'{o}sa-Seymour conjecture, see \cite{KSSz-PS1,KSSz-PS2})
or $H$-factors for a fixed graph $H$ 
(Alon-Yuster conjecture, see \cite{KSSz-AY}).
In this paper we apply this method again.

We will use the following main parameters
\begin{equation}\label{para}
0 < \kappa\ll \eps \ll d \ll \beta \ll \alpha \ll 1,
\end{equation}
where $a\ll b$ means that $a$ is sufficiently small compared to $b$.
For simplicity we do not compute the actual dependencies, although
it could be done. Throughout the rest of the proof we assume that
\begin{equation}\label{k}
2\leq k \leq \kappa n.
\end{equation}

We apply Lemma \ref{degree} for $G$,
with $\eps$ and $d$ as in (\ref{para}).
We get a partition of $V=\cup_{0\leq i\leq l}V_i$.
We define the following so-called {\bf reduced graph} $G_r$:
The vertices of $G_r$ are the clusters $V_i,\,1\leq i \leq l,$
in the partition and there is an edge between two clusters
if they form an $\eps$-regular pair in $G'$
with density exceeding $d$.
Since in $G'$ 
$$\delta(G') > \delta(G) - (d+\eps) n \geq 
\left( \lceil \frac{n}{2} \rceil
+ \lfloor \frac{k}{2} \rfloor  - 1 \right) 
- (d+\eps) n \geq \left( \frac{1}{2} - (d + \eps )\right) n,$$
an easy calculation shows that in $G_r$ we have
\begin{equation}\label{rfok}
\delta(G_r) \geq \left(\frac{1}{2} - 3d \right)l.
\end{equation}
Then Lemma \ref{ep} implies that we can find a matching $M$ in $G_r$
of size at least $\left(\frac{1}{2} - 3d \right)l$. Put 
$|M|=m$. Let us put the vertices of the clusters not covered by $M$
into the exceptional set $V_0$. For simplicity $V_0$ still 
denotes the resulting set. Then 
\begin{equation}\label{kiv}
|V_0| \leq 6 d l L + \eps n \leq 7 d n.
\end{equation}
Denote the $i$-th pair in $M$ by $(V_1^i, V_2^i)$ for $1\leq i \leq m$.
Let $S = v_1, v_2, \ldots , v_k$ be any sequence drawn from $V(G)$.

The rest of the paper is organized as follows. In the next section we
show that if certain extremal conditions hold then we can find directly 
the desired Hamiltonian cycle without using the above method. Then
assuming that the extremal conditions do not hold  we do the following.
In Section 3.3 first we find a short $S$-path $P$. 
Then we find short connecting paths between
the consecutive edges in the matching $M$ (for $i=m$ the next edge is
$i=1$). The first connecting path between $(V_1^1, V_2^1)$ and
$(V_1^2, V_2^2)$ will also contain $P$, the others have length
exactly 3. In Section 3.4 we will take care of the exceptional vertices
and make some adjustments by extending some of the connecting paths so
that the distribution of the remaining vertices inside each
edge in $M$ is perfect, i.e. there are the same number of vertices
left in both clusters of the edge. Finally applying Lemma \ref{ham}
we close the Hamiltonian cycle in each edge. We give the simple
proof of Theorem 2 in Section 4. We finish by some open problems
in Section 5.

\subsection{Extremal cases}

In this section we show that if $G$ satisfies certain extremal
conditions, then we can find the desired Hamiltonian cycle directly.
The first extremal case is the following.

\begin{lemma}\pont\label{extr1}
There exists an $n_3>0$ such that the following holds.
Assume that $G$ is a graph on $n\geq n_3$ vertices satisfying
(\ref{fok}), $k$ satisfies (\ref{k}), and there are $A, B\subset V(G)$ 
such that
\begin{itemize}
\item $A\cap B = \emptyset, |A|, |B| \geq (1 - \alpha ) \frac{n}{2}$,
\item $d(A, B) < \alpha$.
\end{itemize}
Then $G$ is a $k$-ordered Hamiltonian graph.
\end{lemma}

{\bf Proof:} 
First we find the set (denoted by $Exc(A)$) of exceptional vertices $x\in A$
for which $deg(x,B) \geq \sqrt{\alpha} |B|$. The density assumption
implies that the number of these exceptional vertices is at most
$\sqrt{\alpha} |A|$. Similarly we find the set (denoted by $Exc(B)$)
of exceptional vertices $y\in B$ for which $deg(y,A) \geq 
\sqrt{\alpha} |A|$. Again we have $|Exc(B)| \leq \sqrt{\alpha} |B|$. 
We remove the vertices of $Exc(A)$ from $A$ and the vertices of $Exc(B)$
from $B$ and we form
$$E = Exc(A) \cup Exc(B) \cup \left( V(G) \setminus (A\cup B)\right).$$
For each vertex $z\in E$, if $deg(z,A)\geq deg(z,B)$, then we add $z$ to
$A$, and we add $z$ to $B$ in the opposite case. For simplicity we still
denote the resulting sets by $A$ and $B$. It is not hard to see that
in $G|_A$ and in $G|_B$ apart from at most $3 \sqrt{\alpha} n$ 
exceptional vertices all the degrees are at least 
$(1 - \alpha^{1/4} ) \frac{n}{2}$, and the degrees of the exceptional 
vertices are at least $\frac{n}{5}$.

Let $S = v_1, v_2, \ldots ,v_k$ be any sequence drawn from $V(G)$, and let
$$A^* = \{v_1, v_2, \ldots ,v_k\}\cap A \text{and} 
B^* = \{v_1, v_2, \ldots ,v_k\}\cap B.$$
The graph $H$ is the spanning subgraph of $G$ with edges set
$$E(H) = E(G) \cap \left( (B\times (A\setminus A^*)) \cup
(A\times (B\setminus B^*)) \right),$$
that is, the graph obtained from $G$ by removing all edges whose
endpoints are each in $A$ or each in $B$ or each in 
$\{v_1, v_2, \ldots ,v_k\}$. 
We now give a lower bound on the size of an arbitrary vertex cover
$C$ of $E(H)$. By K\"{o}nig's Theorem (see \cite{B})
this gives us the same lower bound on the size of a maximum matching.
If $(A\setminus A^*) \subset C$ or $(B\setminus B^*) \subset C$,
then $|C|\geq k$. Otherwise there exist $x\in (A\setminus A^*) \setminus C$
and $y\in (B\setminus B^*) \setminus C$ and $N(x)\cup N(y) \subset C$.
We have
$$|C|\geq deg_H(x) + deg_H(y) \geq \lceil \frac{n}{2} \rceil + 
\lfloor \frac{k}{2} \rfloor - 1 - (|A| - 1) +
\lceil \frac{n}{2} \rceil + 
\lfloor \frac{k}{2} \rfloor - 1 - (|B| - 1) =$$
$$= 2 \lfloor \frac{k}{2} \rfloor +
2 \lceil \frac{n}{2} \rceil - n.$$
Here the last expression is $k-1$ if $k$ is odd and $n$ is even, and
at least $k$ in all other cases. Therefore, we can always find a matching
$M$ with size $k-1$ if $k$ is odd and with size $k$ if $k$ is even.
We are going to use the edges in $M$ as ``bridges'' on the transitions
between $A$ and $B$.

We are going to use the following fact.

\begin{fact}\pont\label{factegy}
For any $k$ pairs $(u_1,w_1), \ldots ,(u_k,w_k)$ of 
vertices in $A$ (analogously in $B$), there exist $k$ internally
disjoint paths $P(u_1,w_1), \ldots ,P(u_k,w_k)$ of lengths at most 4
such that $P(u_i,w_i)$ has $u_i$ and $w_i$ as endpoints.
\end{fact}

Indeed, this follows from (\ref{k}) and the fact that the number of
exceptional vertices in $A$ is much smaller than the minimum degree
in $G|_A$.

Given our matching $M$ and the sequence $S= v_1, v_2, \ldots ,v_k$, first
we construct an $S$-cycle $C_{k+1}$ (not necessarily Hamiltonian).
Let $v_{k+1} = v_1$. Initially $C_1 = (v_1)$, and the following cases
are applied iteratively for adding $v_{i+1}$ to $C_i$ until $v_{k+1} = v_1$
is added again, in which case $C_{k+1}$ is an $S$-cycle (note that
the intermediate $C_i$-s, $i\leq k$ are only paths):
\begin{itemize}
\item If $\{v_i,v_{i+1}\}\subset A$, then replace $v_i$ by 
$P(v_i,v_{i+1})$ in $C_i$,
yielding $C_{i+1}$.
\item If $v_i\in A$ and $v_{i+1}\in B$ and $v_i$ is incident with an edge
$(v_i,w)$ of $M$, then add $P(w,v_{i+1})$ to $C_i$, yielding $C_{i+1}$.
\item If $v_i\in A$ and $v_{i+1}\in B$ and $v_i$ is not incident with
an edge of $M$, but $v_{i+1}$ is incident with an edge $(w,v_{i+1})$
of $M$, then replace $v_i$ by $P(v_i,w)$ and $v_{i+1}$.
\item If $v_i\in A$ and $v_{i+1}\in B$ and $v_i$ and $v_{i+1}$ are not 
incident with edges of $M$, then there must be an edge $(u,w)\in M$
such that $u\in A\setminus A^*$ and $w\in B\setminus B^*$. Replace
$v_i$ with $P(v_i,u)$ and $P(w,v_{i+1})$ in $C_i$, yielding $C_{i+1}$.
\item The cases for $v_i\in B$ are handled analogously.
\end{itemize}

Here we used the fact that the number of edges in $M$ is at least the number of
transitions we must make between $A$ and $B$, thus we never run out of bridges
in $M$. If $C = C_{k+1}$ contains two consecutive vertices $u_j$ and 
$u_{j+1}$ of $A$ but not all of $A$ belongs to $C$ (an analogous argument
holds for $B$), then let $T$ contain $u_j, u_{j+1}$ and all the vertices 
of $A$ not in $C$. It is not hard to see that the degree conditions in 
Lemma \ref{hamcon} are satisfied in $G|_T$ (with much room to spare)
and thus $G|_T$ is 
Hamiltonian-connected, i.e. there is a Hamiltonian path in $G|_T$ from
$u_j$ to $u_{j+1}$. We add the intermediate vertices of this path
to $C$ between $u_j$ and $u_{j+1}$.

Let $T'$ be the set of vertices not yet in $C$. If $T'=\emptyset$, then we 
are finished. Otherwise either $\{v_1, v_2, \ldots ,v_k\}\subset A$ and
$T'=B$ or $\{v_1, v_2, \ldots ,v_k\}\subset B$ and $T'=A$. Without loss of
generality, assume $T'=A$ and let $x_i$ and $x_j$ two consecutive vertices
in $C$. Choose distinct $y_i, y_j\in B$ such that $(x_i,y_i), (x_j,y_j)\in
E(G)$. Again Lemma \ref{hamcon} implies that $G|_{T'}$ is Hamiltonian-
connected, thus there is a Hamiltonian path in $G|_{T'}$ from $y_i$ to $y_j$.
We add this path between $x_i$ and $x_j$ in $C$ to get the desired
Hamiltonian cycle.

The second extremal case is the following.

\begin{lemma}\pont\label{extr2}
There exists an $n_4>0$ such that the following holds.
Assume that $G$ is a graph on $n\geq n_4$ vertices satisfying
(\ref{fok}), $k$ satisfies (\ref{k}), and there are $A, B\subset V(G)$ 
such that
\begin{itemize}
\item $A\cap B = \emptyset, |A|, |B| \geq (1 - \alpha ) \frac{n}{2}$,
\item $d(A, B) > ( 1 - \alpha)$.
\end{itemize}
Then $G$ is a $k$-ordered Hamiltonian graph.
\end{lemma}

{\bf Proof:}
We start similarly to the proof of Lemma \ref{extr1}.
First we find the set (denoted by $Exc(A)$) of exceptional vertices $x\in A$
for which $deg(x,B) < (1 - \sqrt{\alpha}) |B|$. The density assumption
implies that the number of these exceptional vertices is at most
$\sqrt{\alpha} |A|$. Similarly we find the set (denoted by $Exc(B)$)
of exceptional vertices $y\in B$ for which $deg(y,A) < 
(1 - \sqrt{\alpha}) |A|$. Again we have $|Exc(B)| \leq \sqrt{\alpha} |B|$. 
We remove the vertices of $Exc(A)$ from $A$ and the vertices of $Exc(B)$
from $B$ and we form
$$E = Exc(A) \cup Exc(B) \cup \left( V(G) \setminus (A\cup B)\right).$$
For each vertex $z\in E$, if $deg(z,A)\geq deg(z,B)$, then we add $z$ to
$B$, and we add $z$ to $A$ in the opposite case. For simplicity we still
denote the resulting sets by $A$ and $B$. It is not hard to see that
in $G|_{A\times B}$ apart from at most $3 \sqrt{\alpha} n$ 
exceptional vertices all the degrees are at least 
$(1 - \alpha^{1/4} ) \frac{n}{2}$, and the degrees of the exceptional 
vertices are at least $\frac{n}{5}$.

Without loss of generality, assume that $|A|-|B| = r \geq 0$. We also
know that $r\leq 3 \sqrt{\alpha} n$. Our goal is to achieve $r=0$ since then
we can apply Lemma \ref{bihamcon}. If there is a vertex $x\in A$ for which
\begin{equation}\label{ki}
deg(x,A) \geq \alpha^{1/4} |A|,
\end{equation}
then we remove $x$ from $A$ and add it to $B$. We iterate this procedure
until either there is no more vertex satisfying (\ref{ki}) or $|A|=|B|$.
Let $S=v_1, v_2, \ldots ,v_k$ be any sequence drawn from $V(G)$ and
assume that the first case is true. Since we have 
$\Delta(G|_A) <  \alpha^{1/4} |A|$, (\ref{fok}) and 
Lemma \ref{matching} imply that
$G|_A$ has a $r$-matching $M$ denoted by 
$(u_1,w_1), \ldots ,(u_{r},w_{r})$. 
In case $|A|=|B|$ we have $r=0$ and $M=\emptyset$.
This time the matching $M$ will be used to balance the discrepancy
between $|A|$ and $|B|$.

First we 
construct a short path $P$ that is an $S$-path and contains all the edges 
of $M$. For this purpose we use the following fact (similar to Fact 
\ref{factegy}).

\begin{fact}\pont\label{factketto}
For any $r+k-1$ pairs $(u_1',w_1'), \ldots ,(u_{r+k-1}',w_{r+k-1}')$ 
of vertices, there exist $r+k-1$ internally
disjoint paths $P(u_1',w_1'), \ldots ,P(u_{r+k-1}',w_{r+k-1}')$ 
in $G|_{A\times B}$ of 
lengths at most 5
such that $P(u_i',w_i')$ has $u_i'$ and $w_i'$ as endpoints.
\end{fact}

Then in case $r>0$ the path $P$ is the following (as a sequence of vertices):
$$P = u_1, P(w_1,u_2), P(w_2,u_3), \ldots , P(w_{r-1},u_{r}),$$
$$P(w_{r},v_1), \left(P(v_1,v_2)\setminus v_1\right), \ldots ,
\left(P(v_{k-1},v_k)\setminus v_{k-1}\right).$$
Furthermore, if both $u_1$ and $v_k$ fall into $A$ or into $B$,
then we add one more vertex $v$ from the other set
to the end of $P$ such that $(v_k, v)\in E(G)$.
In case $r=0$ $P$ is the following:
$$P = P(v_1,v_2), \ldots , \left(P(v_{k-1},v_k)\setminus v_{k-1}\right).$$
Again, if both $v_1$ and $v_k$ fall into the same set we add one more vertex
from the other set to the end of $P$.
Lemma \ref{bihamcon} finds the remaining part of the Hamiltonian $S$-cycle
in both cases. Note that for this purpose we could also use Lemma \ref{ham},
the remaining bipartite graph is super-regular with the appropriate choice
of parameters, but here the much simpler Lemma \ref{bihamcon} also
suffices.

An easy consequence of Lemmas \ref{extr1} and \ref{extr2} is our main
extremal case.

\begin{lemma}\pont\label{extr3}
There exists an $n_5>0$ such that the following holds.
Assume that $G$ is a graph on $n\geq n_5$ vertices satisfying
(\ref{fok}), $k$ satisfies (\ref{k}), and there are $A, B\subset V(G)$ 
(not necessarily disjoint) such that
\begin{itemize}
\item $(1 - \beta ) \frac{n}{2} \leq |A|, |B| \leq \frac{n}{2}$,
\item $d(A, B) < \beta$.
\end{itemize}
Then $G$ is a $k$-ordered Hamiltonian graph.
\end{lemma}

{\bf Proof:} We have three cases:

Case 1: $|A\cap B| \geq (1 - \sqrt{\beta} ) \frac{n}{2}$.
In this case the statement follows from Lemma \ref{extr2}.

Case 2: $\sqrt{\beta} \frac{n}{2} \leq 
|A\cap B| < (1 - \sqrt{\beta} ) \frac{n}{2}$. 
This case is not possible under the given assumptions.

Case 3: $|A\cap B| < \sqrt{\beta} \frac{n}{2}$. 
The statement follows from Lemma \ref{extr1}.

\subsection{Connecting paths}

In the remainder of the proof of Theorem 1 we may assume that the extremal 
conditions in Lemma \ref{extr3} do not hold, since otherwise we can find 
the desired Hamiltonian cycle directly. For constructing the
connecting paths we are going to use the following fact several times.
\begin{fact}\pont\label{sokut}
If $x, y \in V(G)$ then there are at least $d n$ internally disjoint
paths of length 3 connecting $x$ and $y$. 
\end{fact}
Indeed, we apply Lemma \ref{extr3} with $A$ and $B$ chosen as follows:
$A$ is an arbitrary subset of $N_G(x)$ with 
$|A|=\lfloor \frac{n}{2} \rfloor$ and $B$ is an arbitrary subset of $N_G(y)$
with $|B|=\lfloor \frac{n}{2} \rfloor$. 

Let $S=v_1, v_2, \ldots ,v_k$ be any sequence drawn from $V(G)$.
Again first we construct a short $S$-path $P$. 
Applying Fact \ref{sokut}, we first connect
$v_1$ and $v_2$ with a path of length 3, then we connect $v_2$ and $v_3$
with a path of length 3 that is internally disjoint from the first connecting
path between $v_1$ and $v_2$, etc. finally we connect $v_{k-1}$ and 
$v_k$ with a path of length 3 that is internally disjoint from all
the connecting paths constructed so far.

For the first connecting path $P_1$ between $(V_1^1, V_2^1)$ and
$(V_1^2, V_2^2)$, first we connect a typical vertex $u$ of $V_2^1$
(more precisely a vertex $u$ with 
$deg(u, V_1^1) \geq ( d - \eps ) |V_1^1|$, most vertices in $V_2^1$
satisfy this) and $v_1$ with a path of length 3, and then we connect
$v_k$ and a typical vertex $w$ of $V_1^2$ (so 
$deg(w, V_2^2) \geq ( d - \eps ) |V_2^2|$) with a path of length 3.
To construct the second connecting path $P_2$ between 
$(V_1^2, V_2^2)$ and
$(V_1^3, V_2^3)$ we just connect a typical vertex of $V_2^2$ and a
typical vertex $V_1^3$ with a path of length 3. Continuing in this fashion,
finally we connect a typical vertex of $V_2^m$ with a typical vertex of
$V_1^1$ with a path of length 3. Thus $P_1$ has length $3(k+1)$, all
other $P_i$-s have length 3. 

We remove the vertices on these connecting
paths
from the clusters, but for simplicity we keep the notation for the
resulting clusters. 
These connecting paths will be parts of the final Hamiltonian cycle.
If the number of remaining vertices
(in the clusters and in $V_0$) is odd, then we take another typical 
vertex $w$ of $V_1^2$ and we extend $P_1$ by a path of length 3
that ends with $w$. So we may always assume that the number of 
remaining vertices is even.

\subsection{Adjustments and the handling of the exceptional vertices}

We already have an exceptional set $V_0$ of vertices in $G$.
We add some more vertices to $V_0$ to achieve super-regularity. From 
$V_1^i$ (and similarly from $V_2^i$)
we remove all vertices $u$ for which
$deg(u, V_2^i) < ( d - \eps ) |V_2^i|$.
$\eps$-regularity guarantees that at most $\eps |V_1^i|\leq \eps L$ such
vertices exist in each cluster $V_1^i$.

Thus using (\ref{kiv}), we still have
$$|V_0| \leq 7 d n + 2\eps n \leq 9 d n.$$
Since we are looking for a Hamiltonian cycle, we have to include the
vertices of $V_0$ on the Hamiltonian cycle as well. We are going to extend
some of the connecting paths $P_i$, so now they are going to contain
the vertices of $V_0$. Let us consider the first vertex 
(in an arbitrary ordering of the vertices in $V_0$) $w$ in $V_0$.
We find a pair $(V_1^i, V_2^i)$ such that either
\begin{equation}\label{egy}
deg(w, V_1^i) \geq d |V_1^i|,
\end{equation}
or
\begin{equation}\label{ketto}
deg(w, V_2^i) \geq d |V_2^i|.
\end{equation}
We assign $w$ to the pair $(V_1^i, V_2^i)$. We extend $P_{i-1}$
(for $i=1$, $P_m$) in $(V_1^i, V_2^i)$ by a path of length 3 in case
(\ref{egy}) holds, and by a path of length 4 in case (\ref{ketto})
holds, so that now the path ends with $w$. To finish the procedure
for $w$, in case (\ref{egy}) holds we add one more vertex $w'$
to $P_{i-1}$ after $w$ such that $(w, w')\in E(G)$ and $w'$
is a typical vertex of $V_1^i$, so
$deg(w', V_2^i) \geq (d - \eps ) |V_2^i|$.
In case (\ref{ketto}) holds we add two more vertices $w', w''$
to $P_{i-1}$ after $w$ such that $(w, w'), (w', w'')\in E(G)$,
$w'$ is a typical vertex of $V_2^i$ and $w''$
is a typical vertex of $V_1^i$.

After handling $w$, we repeat the same procedure for the other
vertices in $V_0$. However, we have to pay attention to several
technical details. First, of course in repeating this procedure
we always consider the remaining free vertices in each cluster;
the vertices on the connecting paths are always removed.
Second, we make sure that we never assign too many vertices of $V_0$
to one pair $(V_1^i, V_2^i)$. It is not hard to see (using 
$d\ll 1$) that we can guarantee that we always assign
at most $\sqrt{d} |V_1^i|$ vertices of $V_0$ to a pair
$(V_1^i, V_2^i)$. Finally, since we are removing vertices from
a pair $(V_1^i, V_2^i)$, we might violate the super-regularity.
Note that we never violate the $\eps$-regularity. Therefore, we do
the following. After handling (say) $\lfloor d^2 n \rfloor$ vertices 
from $V_0$, we
update $V_0$ as follows. In a pair $(V_1^i, V_2^i)$ we remove
all vertices $u$ from $V_1^i$ (and similarly from $V_2^i$)
for which 
$deg(u, V_2^i) < (d - \eps ) |V_2^i|$
(again, we consider only the remaining vertices).
Again, we added at most $2 \eps n$ vertices to $V_0$. In $V_0$
we handle these vertices first and then we move on to the other vertices
in $V_0$.

After we are done with all the vertices in $V_0$, we might have a small
discrepancy ($\leq 2 \sqrt{d} |V_1^i|$) among the remaining vertices
in $V_1^i$ and in $V_2^i$ in a pair. Therefore, we have to make
some adjustments. Let us take a pair $(V_1^i, V_2^i)$ with a
discrepancy $\geq 2$ (if one such pair exists), say $|V_1^i|\geq 
|V_2^i| + 2$ (only remaining vertices are considered). Using Lemma \ref{extr3}
we find an {\bf alternating path} (with respect to $M$) in $G_r$
of length 6 starting with $V_1^i$ and ending with $V_2^i$. 
Let us denote this path by
$$V_1^i, V_2^{i_1}, V_1^{i_1}, V_1^{i_2}, V_2^{i_2}, V_1^i, V_2^i$$
(the construction is similar if the clusters in
$(V_1^{i_1}, V_2^{i_1})$ or in $(V_1^{i_2}, V_2^{i_2})$
are visited in different order). 
We remove a typical vertex from $V_1^i$ and we add it
to $V_1^{i_1}$, then we remove a typical vertex from $V_1^{i_1}$ and we add
it to $V_2^{i_2}$, finally we remove a typical vertex from
$V_2^{i_2}$ and we add it to $V_2^i$. When we add a new vertex to a
pair $(V_1^j, V_2^j)$, we extend the connecting path $P_{j-1}$
by a path of length 4 in the pair so that it now includes the new vertex.

Now we are one step closer to the perfect distribution, and by
iterating this procedure we can assure that the discrepancy in
every pair is $\leq 1$. We consider only those pairs for
which the discrepancy is exactly 1, so in particular the number
of remaining vertices in one such a pair is odd. From the construction
it follows that we have an even number of such pairs. We pair up 
these pairs arbitrarily. If $(V_1^i, V_2^i)$ and $(V_1^j, V_2^j)$
is one such pair with $|V_1^i| = |V_2^i|+1$ and 
$|V_1^j| = |V_2^j|+1$ (otherwise similar), then similar to the
construction above, we find an alternating path in $G_r$ of
length 6 between $V_1^i$ and $V_2^j$, and we move a typical vertex
of $V_1^i$ through the intermediate clusters to $V_2^j$.

Thus we may assume that the distribution is perfect, in every
pair $(V_1^i, V_2^i)$ we have the same number of vertices left.
In this case Lemma \ref{ham} closes the Hamiltonian cycle in every pair.

\section{Proof of Theorem 2}

We consider the graph $G$ with vertices
$$\left\{u_1, \ldots , u_{\lfloor \frac{n}{2} \rfloor}, w_1 , \ldots , 
w_{\lceil \frac{n}{2}
\rceil} \right\}$$ such that 
$U = \left\{u_1, \ldots , u_{ \lfloor \frac{n}{2} \rfloor} \right\}$ and
$W = \left\{w_1 , \ldots , w_{\lceil \frac{n}{2}
\rceil} \right\}$ induce complete subgraphs of $G$. The edges of $G$
between $U$ and $W$ are
$$\left( U\times \left\{w_1 , \ldots , w_{\lfloor \frac{k}{2}
\rfloor} \right\}\right) \cup
\left( W\times \left\{u_1 , \ldots , u_{\lfloor \frac{k}{2} \rfloor - 1} 
\right\}\right).$$
It is easily seen that
$$\delta(G) = \lceil \frac{n}{2} \rceil + \lfloor \frac{k}{2} \rfloor - 2,$$
as required.
Furthermore $G$ does not contain a 
Hamiltonian cycle which encounters
$$u_{\lfloor \frac{k}{2}
\rfloor} - w_{\lfloor \frac{k}{2}
\rfloor + 1} - u_{\lfloor \frac{k}{2}
\rfloor + 1} - w_{\lfloor \frac{k}{2}
\rfloor + 2} - \ldots - u_{2 \lfloor \frac{k}{2}
\rfloor - 1}-
w_{2 \lfloor \frac{k}{2}
\rfloor} \left( - u_{2 \lfloor \frac{k}{2}
\rfloor} \text{if $k$ is odd} \right)$$ in this order.
This follows from the fact, that every $u_i - w_{i+1}$ and $w_i - u_i$
transition uses at least one vertex from the vertices
$$\left\{w_1 , \ldots , w_{\lfloor \frac{k}{2}
\rfloor} \right\} \cup 
\left\{u_1 , \ldots , u_{\lfloor \frac{k}{2} \rfloor - 1} 
\right\}.$$
However, the number of transitions is always more than the number 
of vertices in this set.
Here we also used the fact that we have enough vertices in $U$ and $W$ since 
$k\leq \lfloor \frac{n}{2} \rfloor$.
Thus $G$ is not a 
$k$-ordered Hamiltonian
graph, finishing the proof of Theorem 2.

\section{Open Problems}
The obvious open problem is to determine the best possible minimum
degree requirement for {\em every} $k, n$ pair satisfying $2\leq k \leq n$
and guaranteeing that the graph is a $k$-ordered Hamiltonian graph.
Theorem \ref{tetel1} is valid only for $n\geq n_0$ and $2\leq k \leq 
\kappa n$.

Another open problem is to determine the best possible Ore-type condition.
In \cite{NS} it is shown that if $n\geq 3$, $3\leq k \leq n$ and
$deg(u) + deg(v) \geq n + 2k -6$ for every pair $u, v$ of nonadjacent
vertices of $G$, then $G$ is a $k$-ordered Hamiltonian graph.

\section{Acknowledgement}
We thank Gary Chartrand and Michelle Schultz for their help and encouragement.

\end{document}